\newmathalphabet*{\E}{eus}{m}{n}

\renewcommand{\P}{{\Bbb P}}
\newcommand{\A}{{\Bbb A}}
\newcommand{\Z}{{\Bbb Z}}
\newcommand{\C}{{\Bbb C}}
\newcommand{\Spec}{\operatorname{Spec}}
\renewcommand{\O}{{\E O}}

\newtheorem{thm}{Theorem}
\newtheorem{lem}{Lemma}
\newtheorem{prop}{Proposition}
\newtheorem{cor}{Corollary}
\newtheorem{defn}{Definition}

\newenvironment{pf}{{\bf Proof}}{$\Box$}

\documentstyle[amscd,amssymb,12pt,righttag]{amsart}

\title{Irreducibility of $\overline{M}_{0,n}(G/P,\beta)$}

\author{Jesper Funch Thomsen}

\address
{Matematisk Institut\\
Aarhus Universitet\\
Ny Munkegade\\
DK-8000 \AA rhus C\\
Denmark}
\email{funch@@mi.aau.dk}
\begin{document}

\maketitle

\bigskip

\section{Introduction}
Let $G$ be a complex connected linear algebraic group,
$P$ be a parabolic subgroup of $G$ and $\beta \in A_1(G/P)$
be a 1-cycle class in the Chow group of G/P. 
An $n$-pointed genus $0$ stable map into $G/P$ representing
the class $\beta$, consists of data $( \mu : C \rightarrow X ;
p_1, \dots , p_n )$, where $C$ is a connected, at most nodal,
complex projective curve of arithmetic genus $0$, and $\mu$ 
is a complex morphism such that $\mu_*[C] = \beta$ in $A_1(G/P)$.
In addition $p_i$, $i=1, \dots, n$ denote $n$ nonsingular 
marked points on $C$ such that every component of $C$, which by 
$\mu$ maps to a point, has at least 3 points which is either nodal
or among the marked points (this we will refer to as every component
of $C$ being stable). The set of $n$-pointed genus $0$ stable maps into 
$G/P$ representing the class $\beta$, is parameterized by a coarse 
moduli space $\overline{M}_{0,n}(G/P,\beta)$.
In general it is known that  $\overline{M}_{0,n}(G/P,\beta)$ is a 
normal complex projective scheme with finite quotient singularities.
In this paper we will prove that  $\overline{M}_{0,n}(G/P,\beta)$ is
irreducible. It should also be noted that we in addition 
will prove that the boundary divisors
 in  $\overline{M}_{0,n}(G/P,\beta)$
, usually denoted by 
$D(A,B,\beta_1, \beta_2)$ ($\beta = \beta_1 + \beta_2$, $A \cup B$ 
a partition of $\{ 1, \dots , n \}$)
, are irreducible.

After this work was carried out we learned that B. Kim and 
R. Pandharipande \cite{KimPan} had proven the same results, and even
proved connectedness of the corresponding moduli spaces
in higher genus. 
Our methods however differs in many ways. 
For example in this paper we consider the action of a Borel
subgroup of $G$ on  $\overline{M}_{0,n}(G/P,\beta)$, while 
Kim and Pandharipande mainly concentrate on maximal torus action. 
Another important difference is that we in this presentation 
proceed by induction on $\beta$. This means that 
the question of  $\overline{M}_{0,n}(G/P,\beta)$ being irreducible,
can be reduced to simple cases.

This work was carried out while I took part in the program
``Enumerative geometry and its interaction with theoretical
physics'' at the Mittag-Leffler Institute. I would like to 
use this opportunity to thank the Mittag-Leffler Institute
for creating a stimulating atmosphere.
Thanks are also due to Niels Lauritzen and S{\o}ren Have 
Hansen for useful discussions concerning the generalization
from full flag varieties to patial flag varieties.
  
\section{Summary on  $\overline{M}_{0,n}(G/P,\beta)$}
\label{summary}

In this section we will summarize the properties 
of the coarse moduli space  $\overline{M}_{0,n}(G/P,\beta)$
which we will make use of. The notes on quantum cohomology
by W. Fulton and R. Pandharipande \cite{FulPan} will serve as our main
reference. 

As mentioned in the introduction the moduli 
space  $\overline{M}_{0,n}(G/P,\beta)$ parameterizes $n$-pointed
genus $0$ stable maps into $G/P$ representing the class $\beta$.
By definition $\beta$ is effective if it is represented 
by some $n$-pointed genus $0$ stable map. 
In the following we will only consider values of 
$n$ and $\beta$ where $\overline{M}_{0,n}(G/P,\beta)$ is non-empty.
This means $\beta$ must be effective and $n \geq 0$, and if 
$\beta = 0$ we must have $n \geq 3$.

The moduli space  $\overline{M}_{0,n}(G/P,\beta)$ is known to be a
normal projective scheme (see \cite{FulPan}). This implies that
 $\overline{M}_{0,n}(G/P,\beta)$ splits up into a finite 
disjoint union of its components. This we will use several times.

\subsection{Contraction morphism}
On  $\overline{M}_{0,n+1}(G/P,\beta)$  we have a contraction morphism
$$ \overline{M}_{0,n+1}(G/P,\beta) \rightarrow  \overline{M}_{0,n}(G/P,\beta)$$
which ``forget'' the $(n+1)$'th marked point. The contraction morphisms
value on a closed point in $\overline{M}_{0,n+1}(G/P,\beta)$, represented
by $(\mu : C \rightarrow G/P; p_1, \dots , p_{n+1})$, is the point in 
 $\overline{M}_{0,n}(G/P,\beta)$ represented by 
$( \mu^{\circ} : C^{\circ} \rightarrow  G/P; p_1 , \dots , p_n)$, 
where $C^{\circ}$ denote $C$ 
with the unstable components collapsed, and $\mu^{\circ}$ is the 
map induced from $\mu$.
From the construction of  $\overline{M}_{0,n}(G/P,\beta)$
it follows, that the contraction map is a surjective map with connected
fibres.

\subsection{Evaluation map}
For each element $a \in \{1 , \dots n \} $ we have an evaluation map
$$  \delta_a :\overline{M}_{0,n}(G/P,\beta) \rightarrow G/P.$$
Its value on a closed point in  $\overline{M}_{0,n}(G/P,\beta)$
corresponding to $(\mu :C \rightarrow G/P ; p_1, \dots , p_n)$ is
defined to be $\mu(p_a)$. 

\subsection{Boundary}
By a boundary point in $\overline{M}_{0,n}(G/P,\beta)$ we will mean
a point which correspond to a reducible curve. Let 
$ A \cup B =\{1, \dots n \}$ be a partition of $\{ 1, \dots n \} $
 in  disjoint sets,
and let $\beta_1, \beta_2 \in A_1(X)$ be effective classes such that
 $\beta = \beta_1 + \beta_2$.
We will only consider the cases when $\beta_1 \neq 0$ (resp. 
$\beta_2 \neq 0$) or $|A| \geq 2$ (resp. $|B| \geq 2$).
With these conditions on $\beta_1, \beta_2, A$ and $B$
we let $ D(A,B,\beta_1, \beta_2)$ 
denote the set of elements in $\overline{M}_{0,n}(G/P,\beta)$ where the 
corresponding curve $C$ is of the following form : 

\begin{itemize}
 \item $C$ is the union of (at most nodal) curves $C_A$ and $C_B$ meeting in a 
       point.
 \item The markings of $A$ and $B$ lie on $C_A$ and $C_B$ 
       respectively.
 \item $C_A$ and $C_B$ represent the classes $\beta_1$ and $\beta_2$
       respectively.
\end{itemize}

Notice here that our restrictions on $A$,$B$,$\beta_1$ and $\beta_2$
is the stability conditions on $C_A$ and $C_B$.

\noindent It is clear that every boundary element lies in at least one of 
these $D(A,B,\beta_1,\beta_2)$. 
The sets $D(A,B,\beta_1, \beta_2)$ are in fact closed,
and we will regard them as subschemes of 
$\overline{M}_{0,n}(G/P,\beta)$  by giving them the reduced scheme structure.
Closely related to $D(A,B,\beta_1, \beta_2)$ is the
scheme $M(A,B, \beta_1, \beta_2)$ defined by the fibre 
square

$$ 
\begin{CD}
M(A,B,\beta_1, \beta_2) @>p_2>> 
\overline{M}_{0,A \cup \{ \centerdot \}
}(G/P,\beta_1) \times \overline{M}_{0,B \cup \{ \centerdot \}
}(G/P,\beta_2) \\
@Vp_1VV   @VV\delta_{\centerdot}^A \times \delta_{\centerdot}^BV \\
G/P @>\Delta >> G/P \times G/P
\end {CD}
$$

Here $\Delta$ is the diagonal embedding and $\delta_{\centerdot}^A$
and $\delta_{\centerdot}^B$ denotes the evaluation maps with respect to
the point $ \{ \centerdot \}$. In \cite{FulPan} it is proved that
$M(A,B,\beta_1,\beta_2)$ is a normal projective variety
and that we have a canonical map
$$ M(A,B,\beta_1,\beta_2) \longrightarrow D(A,B,\beta_1, 
\beta_2). $$
This map is clearly surjective.
As
$M(A,B,\beta_1,\beta_2)$ is a closed subscheme of $\overline{M}_{0,A \cup 
\{ \centerdot \} }(G/P,\beta_1) \times  \overline{M}_{0,B \cup 
\{ \centerdot \} }(G/P,\beta_2)$ we can regard the closed points
 of $M(A,B,\beta_1, \beta_2)$ as elements of the form $(z_1,z_2)$,
where $z_1 \in \overline{M}_{0,A \cup \{ \centerdot \} }(G/P,\beta_1)$ and 
$z_2 \in \overline{M}_{0,B \cup \{ \centerdot \} }(G/P,\beta_2)$. The image 
of $(z_1,z_2)$ in $D(A,B,\beta_1,\beta_2)$ will then be denoted 
by $z_1 \sqcup z_2$. Given  $z_1 \in \overline{M}_{0,A \cup \{ \centerdot
 \} }(G/P,\beta_1)$, $z_2 \in \overline{M}_{0,B \cup \{ \centerdot \} \cup
\{ * \}}(G/P,\beta_2)$ and  $z_3 \in \overline{M}_{0,C \cup \{ * \}
 }(G/P,\beta_3)$
, with $\delta_{\centerdot}(z_1) = \delta_{\centerdot}(z_2)$ and 
$\delta_*(z_2) = \delta_*(z_3)$,
 we then have the identity $(z_1 \sqcup z_2)
 \sqcup z_3 = z_1 \sqcup (z_2 \sqcup z_3)$ inside $ \overline{M}_{0,A 
\cup B \cup C}(G/P,\beta_1+ \beta_2 +\beta_3)$.

\subsection{G-action}
As mentioned in the introduction we have a $G$-action
$$G \times  \overline{M}_{0,n}(G/P,\beta) \rightarrow  
\overline{M}_{0,n}(G/P,\beta).$$
On closed points we can describe the action in the following 
way. Let $x \in  \overline{M}_{0,n}(G/P,\beta)$ be a closed point 
corresponding to the data $(\mu : C \rightarrow G/P ; p_1 , \dots , p_n )$,
 and let $g$ be a closed point in $G$. Then $g \cdot x$ is the point 
in  $\overline{M}_{0,n}(G/P,\beta)$ corresponding to 
 $(\mu_g : C \rightarrow G/P ; p_1 , \dots , p_n )$, where $\mu_g = 
(g \cdot) \circ \mu$. Here $g \cdot$ denotes multiplication with
$g$ on $G/P$. 

\subsection{Special cases}
The following special cases of our main result follows from the 
construction and formal properties of our moduli spaces.

\noindent $\beta = 0$ : Here the moduli space 
 $\overline{M}_{0,n}(G/P,\beta)$ is canonical isomorphic to
  $\overline{M}_{0,n} \times G/P$, where $\overline{M}_{0,n}$
denote the moduli space of stable $n$-pointed curves of genus $0$.
As  $\overline{M}_{0,n}$ is known to be irreducible \cite{Knudsen}
we get that 
$\overline{M}_{0,n}(G/P,0)$ is irreducible.

\noindent $G/P = \P^1$ : The irreducibility of $\overline{M}_{0,n}(\P^1,d)$
follows from the construction of the moduli space in \cite{FulPan}.
First of all $\overline{M}_{0,0}(\P^1, 1) \cong \Spec (\C)$ 
so we may assume that $(n,d) \neq (0,1)$. With this assumption 
$\overline{M}_{0,n}(\P^1, d)$ is the quotient of a variety $M$ 
by a finite group. Now $M$ is glued together by the moduli spaces
$\overline{M}_{0,n}(\P^1,d,\overline{t})$ of $\overline{t}$-maps
spaces (here $\overline{t}=(t_0,t_1)$ is a basis of $\O_{\P^1}(1)$).
See section 3 in \cite{FulPan} for a  definition of
 $\overline{M}_{0,n}(\P^1, d, \overline{t})$. The moduli spaces 
 $\overline{M}_{0,n}(\P^1, d, \overline{t})$ are irreducible
(in fact they are $\C^*$-bundles over an open subscheme of
$\overline{M}_{0,m}$ for a suitable $m$). This follows from
the proof of Proposition 3.3 in \cite{FulPan}. 
 It is furthermore
clear that  $\overline{M}_{0,n}(\P^1, d, \overline{t})$, 
and  $\overline{M}_{0,n}(\P^1, d, \overline{t}')$ intersect 
non-trivially for different choices of bases $\overline{t}$ and
$\overline{t}'$. This imply that $\overline{M}_{0,n}(\P^1,d)$ 
is connected, and as it is locally normal
it must be irreducible.

\section{Flag varieties}
\label{flag}

In this section we will give a short review on flag varieties. 
Main references will be \cite{Springer}, \cite{Demazure} and 
\cite{Kock}. In \cite{Springer} one can find the general
theory on the structure of linear algebraic groups. The Chow
group of $G/B$, where $B$ is a Borel subgroup, can be found in 
\cite{Demazure}. From this one easily recovers the Chow group
for a general flag variety $G/P$ (e.g. \cite{Kock} Section 1).

\subsection{Schubert varieties.}
Let $G$ be a complex connected linear algebraic group and 
$P$ be a parabolic subgroup of $G$. As we will only be interested 
in the quotient $G/P$, we may assume that $G$ is semisimple.
Fix a maximal torus $T$ and a Borel subgroup $B$ such that 
$$ T \subseteq B \subseteq P \subseteq G. $$
Let $W$ (resp. $R$) denote the Weyl group (resp. roots) associated to 
$T$ and let $R^+$ denote the positive roots with respect to $B$.
Let further $D \subseteq R^+$ denote the simple 
roots.
Given $\alpha \in R$ we let $s_{\alpha} \in W$ denote the corresponding
reflection.

From general theory on algebraic groups we know that 
$P$ is associated to a unique subset $I \subseteq D$, such that 
$P = B W_I B$, where $W_I$ is the subgroup of $W$ generated by 
the reflections $s_{\alpha}$ with $\alpha \in I$. 
The flag variety $G/P$ is then the disjoint union of a finite number
of $B$-invariant subsets $C(w) = B w P / P$ with $w \in W^I$, where 
$$W^I = \{ w \in W | w \alpha \in R^+ \text{ for all } \alpha \in I \}.$$
Each $C(w)$, $w \in W^I$ is isomorphic to $\A^{l(w)}$. Here $l(w)$
denotes the length of a shortest expression of $w$ as a product 
of simple reflections $s_{\alpha}$, $\alpha \in D$. The closures 
of $C(w)$, $w \in W$, inside
$G/P$ is called the generalized Schubert varieties.
We will denote them by $X_w$, $w \in W$, respectively.
In case $l(w)=1$ we have $X_w \cong \P^{\,1}$.

\subsection{Chow group.}
The Chow group $A_*(G/P)$ is freely generated. As a basis we can pick 
$[X_w]$, $w \in W^I$. In \cite{Kock} it is proved that this basis 
is orthogonal. Using that positive classes intersect in positive 
classes on G/P (Cor. 12.2 in \cite{FulInt}), we conclude that  
a class in  
$A_*(G/P)$ is positive (or zero) if and only if it is of the form 
$$ \sum _{w \in W^I} a_w [X_w] \text{ with } a_w \geq 0.$$

\subsection{Effective classes.}
Let $\beta \in A_1(G/P)$. 
From above it is clear that $\beta$
can only be effective (in the 
sense of  Section \ref{summary}), if $\beta$ is a positive linear 
combination of $[X_{s_{\alpha}}]$ with $\alpha \in D \cap W^I
= D \setminus I$. Noticing that $X_{s_{\alpha}} \cong \P^{\,1}$, 
$\alpha \in D \setminus I$,  implies 
the inverse, that is, a positive linear combination of $[
{X_{s_{\alpha}}}]$, $\alpha \in D \setminus I$ is effective.

Using the above we can introduce a partial ordering on the set of effective
classes in $A_1(G/P)$.

\begin{defn}
Let $\beta_1$ and $\beta_2$ be effective classes. 
If there exist an effective class $\beta_3$ such that 
$\beta_2 = \beta_1 + \beta_3$ we write $\beta_1 \prec 
\beta_2$. If $\beta$ is an effective class with the property 
$$\beta' \prec \beta \Rightarrow \beta' = 0 \text{ or } \beta' 
= \beta$$
we say that $\beta$ is irreducible. An effective class $\beta$
is reducible if it is not irreducible.
\end{defn}

Notice that a non-zero effective class $\beta$ is irreducible if and 
only if $\beta = [X_{s_{\alpha}}]$ for some $\alpha \in
D \setminus I$.

In the proof of the irreducibility of  $\overline{M}_{0,n}(G/P,\beta)$
we will use induction on $\beta$ with respect to this ordering. 
This is possible because given an effective class $\beta \in 
A_1(G/P)$, there is only finitely many other effective classes $\beta'$
with $\beta' \prec \beta$.

\subsection{Summary.}
We are ready to summarize what will be important for us 

\begin{itemize}
 \item The set of effective classes in $A_1(G/P)$  has 
       a $\Z_{\geq 0}$-basis represented by $B$-invariant closed subvarieties
       $X_{s_{\alpha}}$, $\alpha \in D \setminus I$, of 
       $G/P$.
 \item The subsets $X_{s_{\alpha}}$, $\alpha \in D
       \setminus I$ are the only $B$-invariant irreducible 
       1-dimensional closed subsets of $G/P$.
 \item $X_{s_{\alpha}} \cong \P^{\,1}$, $\alpha \in D \setminus I$.
\end{itemize}

\section{Boundary of $\overline{M}_{0,n}(X,\beta)$}
\label{boundary}

In this section we begin the proof of our main result. 
Remember that our convention is that whenever we write 
$\overline{M}_{0,n}(G/P,\beta)$, $D(A,B,\beta_1,\beta_2)$
or $M(A,B,\beta_1,\beta_2)$, we assume that these are well defined 
and non-empty. 
From now on we will assume that $G$, a semisimple linear algebraic
group, and a parabolic subgroup $P$ have been fixed. We let $X$ denote
$G/P$.

We will need to know when $D(A,B,\beta_1,\beta_2)$ 
is irreducible and for this purpose we have the following proposition.

\begin{prop}
\label{divisor}
Suppose that $\overline{M}_{0,A 
\cup \{ \centerdot \} }(X,\beta_1)$ and   $\overline{M}_{0,B 
\cup \{ \centerdot \} }(X,\beta_2)$ are irreducible. Then 
the scheme $M(A,B,\beta_1,\beta_2)$ is also irreducible. In particular 
$D(A,B,\beta_1,\beta_2)$ will be irreducible in this case.

\begin{pf}

As $M(A,B,\beta_1,\beta_2)$ is a normal scheme it splits
up into a disjoint union of irreducible components $C_1,C_2, \dots,
C_l$. Our task is to show that $l=1$. Consider the natural map 
$\pi : G \rightarrow G/P$. Locally (in the Zariski topology)
this map has a section (\cite{Jantzen} p.183) , i.e.
there exists an open cover $ \{U_i \}_{i \in I}$ of $X$
(we assume $U_i \neq \emptyset$)  and 
morphisms $s_i : U_i \rightarrow G$ such that $\pi \circ s_i$ 
is the identity map.
By pulling back the covering  $ \{U_i \}_{i \in I}$ of $X$,
by the evaluation maps $\delta_{\centerdot}^A$ and 
$\delta_{\centerdot}^B$, we get open coverings  
$ \{ V_i^A \}_{i \in I} $
and  $ \{ V_i^B \}_{i \in I} $  of $\overline{M}_{0,A 
\cup \{ \centerdot \} }(X,\beta_1)$ and $\overline{M}_{0,B 
\cup \{ \centerdot \} }(X,\beta_2)$ respectively. 
Finally an open cover $ \{W_i \}_{i \in I}$ of 
$M(A,B,\beta_1,\beta_2)$ is obtained
by setting $W_i = p_1^{-1}(U_i)= p_2^{-1}(V_i^A \times V_i^B)$.
We claim
$$\forall i,j \in I : W_i \cap W_j
\neq \emptyset.$$
To see this consider $U_i,U_j \subseteq
X$. As $X$ is irreducible there exists a closed point 
$x \in U_i \cap U_j$. Using that $G$ acts 
transitively on $X$ we can choose elements $z_1 \in 
\overline{M}_{0,A \cup \{ \centerdot \} }(X,\beta_1)$ and $z_2 \in 
\overline{M}_{0,B \cup \{ \centerdot \} }(X,\beta_2)$, with 
$\delta_{\centerdot}^A(z_1) = \delta_{\centerdot}^B(z_2)=x$.
With these choices it is clear that $(z_1,z_2)$ correspond to a point in 
$W_i \cap W_j$.

Next we want to show that $W_i$ is
irreducible. For this consider the map 
$$
\begin{array}{cccc}
\psi_i : & V_i^A \times V_i^B & \rightarrow &
\overline{M}_{0,A 
\cup \{ \centerdot \} }(X,\beta_1) \times \overline{M}_{0,B \cup 
\{ \centerdot \} }(X,\beta_2) 
\vspace{2ex}\\ 
& (z_1,z_2) & \mapsto & (z_1, ((s_i \circ \delta_{\centerdot}^A)(z_1))
((s_i \circ \delta_{\centerdot}^B)(z_2))^{-1} z_2) 
\end{array}
$$

where we use the group action of $G$ on  $\overline{M}_{0,B 
\cup \{ \centerdot \} }(X,\beta_2)$. By definition $\psi_i$ 
factors through $W_i$. We therefore have an induced map
$$\psi_i^{'} : V_i^A \times V_i^B \rightarrow W_i. $$
Clearly $\psi_i^{'} \circ p_2$ is the identity map. This 
implies that $\psi_i^{'}$ is surjective, and as $V_i^A 
\times V_i^B$ is irreducible, we get that $W_i$ is irreducible.
   
At last we notice  
that as $W_i$ is irreducible it must be contained in one of the 
components $C_1,C_2, \dots ,C_l$ of $M(A,B,\beta_1, \beta_2)$.
On the other hand the $W_i$'s intersect non-trivially so all of
them must be contained in the same component. But 
$ \{ W_i \} _{i \in I }$ was an
open cover of $M(A,B,\beta_1,\beta_2)$. We conclude that $l = 1$,
as desired. Being a surjective image of $M(A,B,\beta_1,\beta_2)$
this implies that $D(A,B,\beta_1,\beta_2)$ is also irreducible.
\end{pf} 
\end{prop}

\section{ Properties of the components of  $\overline{M}_{0,n}
(X,\beta)$} 

In this section we study the behaviour of the components of 
$\overline{M}_{0,n}(X,\beta)$. 
Let $K_1, K_2, \dots, 
K_l$ denote the components of $\overline{M}_{0,n}(X,\beta)$. As 
$ \overline{M}_{0,n}(X,\beta)$ is normal, 
the $K_i$'s are disjoint. Remember that we had a group action 
of $G$ on $\overline{M}_{0,n}(X,\beta)$ which was introduced in Section
\ref{summary}. We claim 

\begin{lem}
\label{inv}
Let $K$ be a component of $\overline{M}_{0,n}(X,\beta)$. Then $K$ 
is invariant under the group action of $G$ on $\overline{M}_{0,n}
(X,\beta)$.

\begin{pf}
Let $\eta : G \times \overline{M}_{0,n}(X,\beta) \rightarrow 
\overline{M}_{0,n}(X,\beta)$ denote the group action, and consider
the image $\eta(G \times K)$ of $G \times K$. As $G \times K$ is
irreducible $\eta(G \times K)$ will also be irreducible. This 
means that $\eta(G \times K)$ is contained in a component,  
say $K_1$, of $\overline{M}_{0,n}(X,\beta)$. On the other hand 
$\eta( \{ e \} \times K) \subseteq K$ (here $e$ denotes the identity
 element in $G$) so we conclude that $K=K_1$.
\end {pf}
\end{lem}

The next lemma concerns the boundary of components of $\overline{M}_{0,0}
(X,\beta)$ when $\beta$ is reducible.

\begin{lem}
\label{bound}
Let $\beta$ be a reducible effective class and 
$K$ be a component of $\overline{M}_{0,0}(X,\beta)$.
Then there exist boundary elements in $K$.

\begin{pf}
Assume $K$ do not have boundary elements. Then by definition
of boundary points, each element in $K$ would correspond to an 
irreducible curve. Using Lemma \ref{inv} we have an induced 
$B$-action on $K$. As $K$ is projective, and $B$ is a connected
 solvable linear algebraic group, we can use Borel's fixed point
Theorem (see \cite{Springer} p.159) to conclude that this action
has a fixed point. This means that there exist $z \in K$ such that 
$b z = z$, for all $b \in B$. Let $\P^{\,1} \stackrel{\mu}{\rightarrow}
X$ be the stable curve with its morphism to $X$ which 
correspond to $z$. By definition of the group action of $B$ on
$z$ we conclude that $\mu(\P^{\,1})$ must be a $B$-invariant subset of $X$.
On the other hand $\mu (\P^{\,1})$ is closed, irreducible and of 
dimension 1. By the properties stated in Section \ref{flag}
$\mu(\P^{\,1})$ must be equal to a 1-dimensional 
Schubert variety $X_{s_{\alpha}}$ of $X$.
From this we conclude that 
$\mu_*[\P^{\,1}] = m [X_{s_{\alpha}}]$, where $m$ 
is a positive integer.
As $\beta$ is reducible $m \geqq 2$.
The closed embedding $i : X_{s_{\alpha}}
\rightarrow X$ induces a map
$i_* : \overline{M}_{0,0}(X_{s_{\alpha}}, m
[X_{s_{\alpha}}]) 
\rightarrow \overline{M}_{0,0}(X,\beta)$, where an element $(C \stackrel{f}
{\rightarrow} X_{s_{\alpha}}) \in  \overline{M}_{0,0}
(X_{s_{\alpha}},
m [X_{s_{\alpha}}])$ goes to $i_*(C \stackrel{f}
{\rightarrow} X_{s_{\alpha}}) =
(C\stackrel{i \circ f}{\rightarrow}
X)$. As $X_{s_{\alpha}}$ is isomorphic to $\P^{\,1}$ we know that 
$ \overline{M}_{0,0}(X_{s_{\alpha}}, m 
[{X}_{s_{\alpha}}])$ is irreducible.
On the other hand $z$ is in the 
image of $i_*$ so 
we conclude that $i_*( \overline{M}_{0,0}(X_{s_{\alpha}}, m
[X_{s_{\alpha}}])) \subseteq K$.
But a boundary element in 
$\overline{M}_{0,0}(X_{s_{\alpha}}, m [
X_{s_{\alpha}}])$ is easy to construct by hand
(as $m \geq 2$), 
which gives us the desired contradiction.
\end{pf}
\end{lem}

The following will also be useful.

\begin{lem}
\label{hit}
Let $\beta \in A_1(X)$ be a reducible effective
class and suppose that $\overline{M}_{0,0}(X,\beta')$ is 
irreducible for $\beta' \prec \beta$. 
Furthermore let $K$ be a component of $\overline{M}_{0,0}(X,\beta)$.
Then there exists a 
non-zero irreducible class $\beta'$, with $\beta -\beta'$ effective,
such that $D(\emptyset, \emptyset,
\beta', \beta - \beta') \cap K  \neq \emptyset$.

\begin{pf}
By Lemma \ref{bound} we can choose a boundary point $z \in K$. 
There exists effective classes $\beta_1$ and $\beta_2$ such
that $z \in D(\emptyset, \emptyset, \beta_1, \beta_2)$.
We may assume that $\beta_1$ is reducible. Choose an effective
non-zero irreducible  class  $\beta'$ and an effective class $\beta''$
such that $\beta_1 = \beta' + \beta''$. Choose also 
$z_1 \in \overline{M}_{0, \{ Q_1 \} }(X,\beta')$,
$z_2 \in \overline{M}_{0, \{ Q_1 \} \cup  \{ Q_2 \} }(X,\beta'')$ and  
$z_3 \in \overline{M}_{0, \{ Q_2 \} }(X,\beta_2)$,
such that $\delta_{Q_1}(z_1) = \delta_{Q_1}(z_2)$ and 
$\delta_{Q_2}(z_2) = \delta_{Q_2}(z_3)$. Then
$z_1 \sqcup z_2 \in \overline{M}_{0, \{ Q_2 \} }(X,\beta_1)$
from which we conclude 
$(z_1 \sqcup z_2) \sqcup z_3 \in D(\emptyset, \emptyset,
\beta_1, \beta_2)$. On the other hand 
$z_2 \sqcup z_3 \in \overline{M}_{0, \{ Q_1 \} }(X,\beta - \beta')$
by which we conclude 
$z_1 \sqcup (z_2 \sqcup z_3) \in D(\emptyset, \emptyset,
\beta', \beta - \beta')$.
Using Proposition \ref{divisor} we know that $ D(\emptyset, \emptyset,
\beta_1, \beta_2)$ is irreducible and as 
$z \in  D(\emptyset, \emptyset,\beta_1, \beta_2) \cap 
K$, we must have $ D(\emptyset, \emptyset,\beta_1, \beta_2) 
\subseteq K$, in particular $(z_1 \sqcup z_2) \sqcup z_3
\in K$. On the other hand 
$$ (z_1 \sqcup z_2) \sqcup z_3 = z_1 \sqcup ( z_2 \sqcup z_3)  
\in D(\emptyset, \emptyset,
\beta', \beta - \beta').$$
This proves the lemma.
\end{pf}
\end{lem}

\section{Irreducibility of $\overline{M}_{0,n}(X,\beta)$}

In this section we will prove that the moduli spaces 
$\overline{M}_{0,n}(X,\beta)$ are irreducible. First we notice
that for $\beta \neq 0$ we can restrict our attention to a fixed
$n$.

\begin{lem}
\label{n-inv}
Let $n_1,n_2 \geqq 0$ be integers and $\beta \in A_1(X) \setminus
 \{ 0 \}$ be an effective class. Then $\overline{M}_{0,n_1}(X,\beta)$
is irreducible if and only if  $\overline{M}_{0,n_2}(X,\beta)$ is
irreducible.

\begin{pf}
It is enough to consider the case $n_2=n+1$ and $n_1=n$ for 
a positive integer $n$. The contraction
morphism $f : \overline{M}_{0,n+1}(X,\beta) \rightarrow
\overline{M}_{0,n}(X,\beta)$ which forgets the $(n+1)$'th point
is a surjective map with connected fibres. Let $K_1,K_2, \dots,
K_s$ (resp. $C_1,C_2, \dots, C_t$) be the components of 
 $\overline{M}_{0,n+1}(X,\beta)$ (resp. $\overline{M}_{0,n}(X,\beta)$).
As the components are mutually disjoint and $f$ is surjective
we must have $s \geqq t$. Let us now restrict our attention to one 
of the components of  $\overline{M}_{0,n}(X,\beta)$, say $C_1$. Assume
that $K_1,K_2, \dots, K_r$ ($r \leqq s$) are the components 
which by $f$ maps to $C_1$. It will be enough to show that $r=1$.
Assume $r \geqq 2$. As

$$ C_1 = \bigcup_{i=1}^{r} f(K_i) $$   
and as $C_1$ is irreducible, at least one of the components $K_1
,K_2, \dots, K_r$ maps surjectively onto $C_1$. 
So there must exist a point $x$ in $C_1$ which is in
the image of at least 2 of the components in  $\overline{M}_{0,n+1}(X,\beta)$.
But then the fibre of $f$ over $x$ is not connected, which 
is a contradiction.
\end{pf}
\end{lem}

The idea in proving the irreducibility of  $\overline{M}_{0,n}(X,\beta)$
is to use induction on the class $\beta \in A_1(X)$. By 
this we mean that we will prove that  $\overline{M}_{0,n}(X,\beta)$
is irreducible assuming the same condition is true for $\beta' \prec
\beta$. 
The first step in the induction procedure will be to show 
that  $\overline{M}_{0,0}(X,\beta)$ is irreducible, when $\beta$ is
a non-zero irreducible class.

\begin{lem}
\label{nonred}
Let $\beta$ be a non-zero irreducible class.
Then $\overline{M}_{0,0}(X,\beta)$ is irreducible.

\begin{pf}
As $\beta$ is a non-zero irreducible class,  $\beta$ must be the class of a
 1-dimensional Schubert variety $X_{s_{\alpha}}$. 
Let $K$ be a component 
of  $\overline{M}_{0,0}(X,\beta)$. As in the proof of Lemma \ref{bound}
we have 
a $B$-action on $K$, which by Borel's fixed point theorem 
is forced to have a fixed point. Let $x \in K$ be a fixed point.
As $\beta$ is irreducible $x$ must correspond to an irreducible 
curve, i.e. $x$ correspond to a map of the form $ \P^{\,1} 
\stackrel{\mu}{\rightarrow} X$. The image $\mu(\P^{\,1})$ is  
a closed 1-dimensional $B$-invariant irreducible
subset of $X$. As by  assumption  $\mu_*[\P^{\,1}] = [\overline
{X}_{s_{\alpha}}]$, we conclude that 
$\mu(\P^{\,1}) = X_{s_{\alpha}}$.
Now $X_{s_{\alpha}}$ is isomorphic to $\P^{\,1}$, 
so $\mu$ must be an isomorphism
onto its image. But clearly every map  $ \P^{\,1} 
\stackrel{f}{\rightarrow} X$ with $f(\P^{\,1})= X_{s_{\alpha}}$, 
which is an 
isomorphism onto its image, represent the same point in 
$\overline{M}_{0,0}(X,\beta)$. Above we have shown that this point 
belongs to every component of  $\overline{M}_{0,0}(X,\beta)$. 
Using that the components of  $\overline{M}_{0,0}(X,\beta)$ are disjoint
the lemma follows.
\end{pf}
\end{lem}

Now we are ready for the general case.

\begin{thm}
\label{thm}
Let $\beta \in A_1(X)$ be an effective class and $X = G/P$ be a 
flag variety. Then  $\overline{M}_{0,n}(X,\beta)$ is irreducible
for every positive integer n.

\begin{pf}
The case $\beta = 0 $ is trivial as noted in Section \ref{summary}.
By Lemma \ref{n-inv} we may therefore assume that $n = 0$. As
remarked above we will proceed by induction. Assume
that the theorem has been proven for $\beta '$ with
$\beta ' \prec \beta$. Referring to Lemma \ref{nonred} we may assume
that $\beta$ is reducible. Write $\beta = \sum_{i=1}^{m}
\beta_i$ as a sum of non-zero irreducible effective classes
$\beta_i$. 
Then $m \geq 2$.
We divide into 2 cases.

Assume first that 
$m=2$. So $\beta = \beta' + \beta''$,
where $\beta'$ and $\beta''$ are effective irreducible
classes. In this case every boundary element lie in 
$D(\emptyset , \emptyset, \beta', \beta'')$, which
we by induction know is irreducible (Proposition \ref{divisor}). 
On the other hand do every component of $\overline{M}_{0,0}
(X,\beta)$ contain a boundary point (by Lemma \ref{bound}).
Using that the components of   $\overline{M}_{0,0}(X,\beta)$
are disjoint, the theorem follows in this case.
 
Assume therefore that $m \geq 3$.
For each $i=1, \dots,m$ choose $z_i \in 
\overline{M}_{0,\{ Q_i \} }(X,\beta_i)$, a point such that $\delta_{Q_i}
(z_i)= eP$, where $\delta_{Q_i}$ is the evaluation map onto
$X$. 
Let $Q = \{ Q_1, Q_2, \dots, Q_m\}$, and choose a
point $z_0 \in \overline{M}_{0,Q}(X,0)$ corresponding to 
a curve $C \cong \P^{\,1}$ and a map $\mu : C \rightarrow X$ such that 
$\mu(C) = eP$. Define
$$ z = z_0 \sqcup (\sqcup_{i=1}^{m} 
 z_i) \in \overline{M}_{0,0}(X,\beta).$$
Then  clearly $z \in D(\emptyset , \emptyset, \beta_i, \beta 
-\beta_i)$ for all $i$. Let $K$
be the component of $\overline{M}_{0,0}(X,\beta)$ which
contains $z$. By the induction hypothesis and Proposition 1,  
$D(\emptyset , \emptyset, \beta_i, \beta - \beta_i)$ is
irreducible for all $i$, which implies $D(\emptyset , \emptyset,
\beta_i, \beta -\beta_i) \subseteq K$ for all $i$.
On the other hand, by Lemma \ref{hit}, every component of 
$\overline{M}_{0,0}(X,\beta)$ will intersect at least 
one of the sets $D(\emptyset , \emptyset, \beta_i, \beta
- \beta_i)$. Using, and now for the last time, that the 
components of $\overline{M}_{0,0}(X,\beta)$ are 
disjoint, the theorem follows.
\end{pf}
\end{thm}

\begin{cor}
Let $X=G/P$ be a flag variety. Then the boundary divisors 
$D(A,B,\beta_1,\beta_2)$ of $\overline{M}_{0,n}(X,\beta)$
are irreducible.
\begin{pf}
Use Proposition \ref{divisor} and Theorem \ref{thm}.
\end{pf} 
\end{cor}

\bibliographystyle{amsplain}
\ifx\undefined\bysame
\newcommand{\bysame}{\leavevmode\hbox to3em{\hrulefill}\,}
\fi

\end{document}